\magnification=1200 
\hsize 158truemm
\overfullrule 0cm

\def\ref#1{\lbrack {#1}\rbrack}

\def\ekv#1#2{$${#2}\eqno(#1)$$}
\def\eekv#1#2#3{$$\eqalignno{&{#2}&({#1})\cr &{#3}\cr}$$}

\font\liten=cmr10 at 8pt
\font\stor=cmr10 at 12pt

\def\ably{arbitrarily}

\def\bdd{bounded}

\def\coef{coefficient}
\def\coeff{coefficient}
\def\const{constant}

\def\canform{canonical transformation}
\def\coef{coefficient}
\def\coeff{coefficient}

\def\ctf{canonical transformation}

\def\dop{differential operator}
\def\ef{eigenfunction}
\def\ev{eigenvalue}

\def\eq{equation}

\def\fu{function}

\def\fop{Fourier integral operator}

\def\hol{holomorphic}
\def\hm{homogeneous}
\def\indep{independent}

\def\mfld{manifold}
\def\ml{microlocal}
\def\neigh{neighborhood}
\def\nondeg{non-degenerate}
\def\op{operator}

\def\pb{problem}

\def\per{periodic}
\def\pert{perturbation}
\def\prop{proposition}
\def\Prop{Proposition}
\def\pol{polynomial}
\def\pop{pseudodifferential operator}

\def\res{resonance}
\def\rhs{right hand side}
\def\sa{selfadjoint}
\def\sc{semiclassical}

\def\st{strictly}
\def\stpsh{\st{} plurisubharmonic}

\def\sufly{sufficiently}
\def\tf{transformation}
\def\Th{Theorem}

\def\tf{transform}

\def\ufly{uniformly}

\def\wrt{with respect to}

\def\Re{{\rm Re\,}}
\def\Im{{\rm Im\,}}

\centerline{\stor Resonances associated to a closed hyperbolic trajectory
in dimension 2.}
\medskip
\centerline{\bf Johannes Sj\"ostrand\footnote{*}{\liten Centre de
Math\'ematiques, Ecole Polytechnique, FR-91128 Palaiseau cedex, and URM
7640 de CNRS}}\footnote{}{\liten Key words: \res{}, scattering pole,
hyperbolic, two.}\footnote{}{\liten MSC2000: 32A99, 35P20, 35S99, 81U99,
37Dxx}
\bigskip
\par\noindent \bf Abstract. \liten We consider \res{}s in the
semi-classical limit ($h\to 0$), generated by a single closed hyperbolic
orbit, for an \op{} on ${\bf R}^2$. We determine all such \res{}s in an
$h$-\indep{} domain. As an application we determine all \res{}s 
generated by a saddle point in a fixed disc around the critical energy.
\medskip
\par\noindent \bf R\'esum\'e. \liten Nous consid\'erons les r\'esonances
engendr\'ees par une trajectoire ferm\'ee hyperbolique pour un
op\'erateur sur ${\bf R}^2$ dans la limite semi-classique ($h\to 0$).
Nous d\'eterminons toutes les r\'esonances dans un domaine ind\'ependant
de $h$. Comme une application nous obtenons toutes les r\'esonances
engendr\'ees par un point selle, dans un disque fix\'e autour de
l'\'en\'ergie critique.
\rm
\bigskip
\centerline{\bf 0. Introduction.}
\medskip
This paper is a continuation of [MeSj], where we made the observation
that for a large class of non-\sa{} \sc{} \op{}s in dimension 2, one can
describe the whole spectrum in some $h$-\indep{} domain in the complex
plane. In that paper we also showed how to determine (in 2 dimensions)
all \res{}s (scattering poles) for $-h^2\Delta +V(x)$ generated by a saddle
point of the potential with critical value $E_0$, in a disc of fixed
radius around $E_0$, with small sectors around $E_0+]0,+\infty [$ and
$E_0-i]0,\infty [$ removed. That is a partial improvement of the results of
Kaidi-Kerdelhu\'e [KaKe], whose result in this case gives all \res{}s in
a disc $D(E_0,h^\delta )$ of center $E_0$ and radius $h^\delta $. Here
$\delta >0$ can be any fixed \const{} and $h>0$ is small enough
depending on $\delta $. 

\par In [GeSj] we studied the \res{}s generated by a closed hyperbolic
trajectory at energy $E_1\in{\bf R}$ and determined all such \res{}s in a
rectangle $]E_1-\epsilon _0,E_1+\epsilon _0[-i[0,Ch[$, where $\epsilon
_0>0$ is \sufly{} small and $C>0$ is an \ably{} large \const{}. It has
long been believed that in analogy with [KaKe], who used quantum Birkhoff
normal forms (qBnf) to improve the result of [Sj] (who got the \res{}s in a
disc of radius $Ch$), the result of [GeSj] should have an improvement
giving all \res{}s in $]E_1-\epsilon _0,E_1+\epsilon _0[-i[0,h^\delta [$
for every fixed $\delta >0$. Either one should apply qBnf directly, or
one should apply the qBnf for \fop{}s by [Ia] to the quantum monodromy
\op{} ([SjZw]). (See also [IaSj].) In this paper, we restrict the attention
to the 2-dimensional case and show how to adapt the result of [MeSj] to
get all \res{}s in
$]E_1-\epsilon _0,E_1+\epsilon _0[-i[h^\delta ,\epsilon _1[$ for
$\epsilon _1>0$ \sufly{} small and with $\delta $ equal to any fixed
constant in $]0,1[$. We also develop the necessary qBnf near the closed
trajectory, to get all \res{}s in $]E_1-\epsilon _0,E_1+\epsilon
_0[-i[0,h^\delta [$. (Presumably the latter result could be extended to
all dimensions if we impose the appropriate non-\res{} condition.) Thus
we get all resonances in a small but $h$-\indep{} rectangle
$]E_1-\epsilon _0,E_1+\epsilon _0[-i[0,\epsilon _1[$. See \Th{} 2.3.

\par In section 3 we reexamine the resonances generated by a saddle point
and improve the result from [MeSj] by treating the "missing sectors".
Thus we get all \res{}s in a fixed disc around the critical value.

\par As in [MeSj] and in earlier works on \res{}s ([HeSj], [Sj2]) a basic
ingredient in the proofs is the use of FBI-\tf{}s and corresponding
weighted spaces of \hol{} \fu{}s (in the spririt of [Sj3]). This has now
become a routine and in order to avoid an excessive length of the paper
we have chosen not to review that material here. This also concerns the
setup of the global Grushin problem which is essentially identical to
that of [MeSj]. See also the first of a series of preprints in
preparation with M. Hitrik.
\bigskip

\centerline{\bf 1. Normal forms.}
\medskip
\par We place ourselves in the general frame work of [HeSj] and let 
\ekv{1.1}
{P=\sum_{\vert \alpha \vert \le m}a_\alpha (x;h)(hD_x)^\alpha ,\ x\in{\bf
R}^2}
be a \sc{} formally \sa{} \dop{} satisfying the general assumptions of
[HeSj] which allow to define the \res{}s in some fixed \neigh{} of
$E_0\in{\bf R}$. 

\par We will also assume that 
\ekv{1.2}{x\mapsto a_\alpha (x;h)\hbox{ extend \hol{}ally to a \neigh{} of
}{\bf R}^n}
with 
\ekv{1.3}
{a_\alpha (x;h)\sim\sum_{j=0}^\infty  a_{\alpha ,j}(x)h^j}
\ufly{} for $x$ in any compact subset of the \neigh{} in (1.2). The \sc{}
principal symbol of $P$ is then given by
\ekv{1.4}{p(x,\xi )=\sum_{\vert \alpha \vert \le m}a_{\alpha ,0}(x)\xi
^\alpha .}
\par The standard example we have in mind is 
\ekv{1.5}
{P=-h^2\Delta +V(x),}
where $V\in C^\infty ({\bf R}^n;{\bf R})$, and $V$ extends to a \hol{}
\fu{} in a set $\vert \Im x\vert <C^{-1}\langle x\rangle $, with $\langle
x\rangle =(1+\vert x\vert ^2)^{1/2}$, and tends to 0 when $x\to \infty $ in
that set. Then
$p(x,\xi )=\xi ^2+V(x)$.

\par Recall from [GeSj] that the trapped set $K(E)$, for $E\in{\rm
neigh\,}(E_0,{\bf R})$ is 
\ekv{1.6}
{K(E)=\{ \rho \in p^{-1}(E);\, \exp tH_p(\rho )\not\to \infty ,\,\, t\to
\pm \infty \}.}
Here $H_p={\partial p\over \partial \xi }\cdot {\partial \over \partial
x}-{\partial p\over \partial x}\cdot {\partial \over \partial \xi }$ is the
Hamilton field of $p(x,\xi )$. We also recall that $K(E)$ is contained in
some fixed compact set when $E$ varies in a \neigh{} of $E_0$. 

\par We assume:
\ekv{1.7}
{K(E_0)\hbox{ is (the image of) a simple closed } H_p\hbox{-trajectory
}\gamma (E_0)\hbox{ of period $T(E_0)$.}}
\ekv{1.8}
{\gamma (E_0)\hbox{ is of hyperbolic type}.}
Recall that (1.8) means that the linearized Poincar\'e map  has the
\ev{}s $\lambda (E_0)$, $1/\lambda (E_0)$ with $\lambda (E_0)\in{\bf R}$,
$\vert \lambda (E_0)\vert >1$. In particular, the Poincar\'e map is
non-degenerate (1 is not an \ev{}), so for $E\in{\rm neigh\,}(E_0,{\bf
R})$, we have a closed $H_p$-trajectory $\gamma (E)$ close to $\gamma
(E_0)$ which depends analytically on $E$. For some small fixed $\epsilon
_0>0$, let 
\ekv{1.9}
{\Gamma =\bigcup_{\vert E-E_0\vert <\epsilon _0}\gamma (E).}
Then $\Gamma $ is a symplectic \mfld{} of dimension 2. Let $\Gamma
_+(E),\Gamma _-(E)\subset p^{-1}(E)$ be the unstable (outgoing) and
stable (incoming) \mfld{}s for the $H_p$-flow. We know that they are
hypersurfaces in $p^{-1}(E)$ which intersect transversally along $\gamma
(E)$. They are also Lagrangian \mfld{}s, and 
$$\Gamma _\pm :=\bigcup_{\vert E-E_0\vert <\epsilon _0}\Gamma _{\pm}(E)$$
are involutive \mfld{}s and can be viewed as the outgoing ($+$) and
incoming ($-$) \mfld{}s for the $H_p$-flow in $p^{-1}(]E_0-\epsilon
_0,E_0+\epsilon _0[)$. They intersect transversally along $\Gamma $.

\par If $\Gamma _+$ is orientable ($\lambda _+(E_0)>1$), then so is
$\Gamma _-$ and we can find an analytic real-valued function $\xi $ in a
\neigh{} of $\Gamma $ with
\ekv{1.10}{\xi =0,\,\, d\xi \ne 0,\hbox{ on }\Gamma _+.} Define the \fu{}
$x$ in a neigh{} of $\Gamma $ to be the solution of 
\ekv{1.11}
{H_{\xi }x=1,\ {x_\vert}_{\Gamma _-}=0.}

\par If $\Gamma _+$ is not orientable ($\lambda _+(E_0)<-1$), we can
still define $\xi $ satisfying (1.10), no more single-valued but
double-valued, with the property that ${\rm ext\,}(\xi )=-\xi $, where
${\rm ext\,}(\xi )$ denotes the extension of $\xi $ obtained making one
tour in a \neigh{} of $\gamma (E_0)$ in the forward direction. We then
define $x$ by (1.11) and observe that ${\rm ext\,}(x)=-x$. Notice that
$x\xi $ is a single-valued function. 

\par Choose symplectic coordinates $(t,\tau )$ on $\Gamma $ with $t$
multivalued:
\ekv{1.12}
{{\rm ext\,}(t)=t+2\pi ,}
and $\tau =\tau (E)$ depending only on $p(\rho )=E$. We may assume that
$\tau (E_0)=0$.  Extend $(t,\tau )$ to a \neigh{} of $\Gamma $, by
solving 
\ekv{1.13}
{H_xt=H_\xi t=0,\ H_x\tau =H_\xi \tau =0.}
This is possible since $[H_x,H_\xi ]=H_{\{ x,\xi \} }=0$ (where $\{ x,\xi
\} =H_x\xi $ is the Poisson bracket) and
$H_x,H_\xi
$ span a plane transversal to $\Gamma $ at each point of $\Gamma $. $\tau $
will be constant on each $\Gamma _{\pm}(E)$. 

\par Then $(t,\tau ;x,\xi )$
are symplectic coordinates and since $H_p$ is tangential to $\Gamma _\pm
(E)$, we see that 
$$\partial _tp(t,\tau ,x,0)=\partial _tp(t,\tau ,0,\xi )=0,\ \partial
_xp(t,\tau ,x,0)=0,\ \partial _\xi p(t,\tau ,0,\xi )=0.$$
It follows that 
$$p(t,\tau ,x,0)=p(t,\tau ,0,\xi )=f(\tau ),$$
so that
$$p(t,\tau ,x,\xi )=f(\tau )+\widetilde{\mu }(t,\tau ,x,\xi )x\xi ,$$
which gives
\ekv{1.14}{p(t,\tau ,x,\xi )=f(\tau )+\mu (t,\tau )x\xi +{\cal O}((x,\xi
)^3).}

\par Let $G=\lambda (t,\tau )x\xi $. Then in the sense of formal Taylor
expansions in $x,\xi $, we get
$$p\circ \exp H_G =\sum {1\over k!}H_G^k p=p+H_Gp+{1\over 2}H_G^2p +{\cal
O}((x,\xi )^3).$$
Here 
$$\eqalign{H_G&=x\xi ((\partial _\tau \lambda )\partial _t-(\partial
_t\lambda )\partial _\tau )+\lambda (t,\tau )(x\partial _x-\xi \partial
_\xi ),\cr
H_Gp&=-(\partial _t\lambda )f'(\tau )x\xi +{\cal O}((x,\xi )^3),\cr
H_G^2p&={\cal O}((x,\xi )^3).}$$
Hence 
$$p\circ \exp H_G=f(\tau )+(\mu (t,\tau )-f'(\tau )\partial _t\lambda
(t,\tau ))x\xi +{\cal O}((x,\xi )^3).$$
Choosing $\lambda $ suitably, we get 
$$\mu (t,\tau )-f'(\tau )\partial _t\lambda (t,\tau )=\langle \mu (\cdot
,\tau )\rangle :={1\over 2\pi }\int_0^{2\pi }\mu (t,\tau ) dt.$$
Replacing $p$ by $p\circ \exp H_G$ (which amounts to expressing $p$ in new
symplectic coordinates), we may assume that
\ekv{1.15}
{p(t,\tau ,x,\xi )=f(\tau )+\mu (\tau )x\xi +{\cal O}((x,\xi )^3).}

\par We have 
\ekv{1.16}
{
f'(\tau )=2\pi /T(f(\tau )) ,
}
where $T(E)$ is the period of $\gamma (E)$. In the orientable case, the
linearized Poincar\'e map at energy $E$ has the \ev{}s $\lambda (E)$,
$1/\lambda (E)$, with $\lambda (E)>1$ given by 
\ekv{1.17}
{\vert \lambda (E)\vert =e^{T(E)\mu (\tau )},\ f(\tau )=E.}
The same relation holds for $\lambda (E)<-1$ in the non-orientable case.

\par We next improve (1.15) to \ably{} high order in $(x,\xi )$. Write 
\ekv{1.18}
{p=p_0+p_2+p_3+...,}
where $p_j$ is \hm{} of degree $j$ in $(x,\xi )$, $p_0=f(\tau )$,
$p_2=\mu (\tau )x\xi $. Let $G_j=G_j(t,\tau ,x,\xi )$ be a \hm{} \pol{}
of degree $j\ge 3$ in $(x,\xi )$ with analytic \coeff{}s, and consider
$p\circ
\exp H_{G_j}=\sum_0^\infty {1\over k!}H_{G_j}^kp$. Here,
\ekv{1.19}
{
H_{G_j}p=-(f'(\tau )\partial _t+\mu (\tau )(x\partial _x-\xi \partial
_\xi ))G_j+{\cal O}((x,\xi )^{j+1}), }
where the first term is \hm{} of degree $j$ in $(x,\xi )$, and 
$$H_{G_j}^kp={\cal O}((x,\xi )^{j+(k-1)(j-2)}),$$
for $k\ge 2$. Consequently,
$$p\circ \exp H_{G_j}=p-(f'(\tau )\partial _t+\mu (\tau )(x\partial _x-\xi
\partial _\xi ))G_j+{\cal O}((x,\xi )^{j+1}).$$

\par In the orientable case, the \eq{}
$$(f'(\tau )\partial _t+\mu (\tau )(x\partial _x-\xi \partial _\xi
))u(t,\tau )x^\alpha \xi ^\beta =v(t,\tau )x^\alpha \xi ^\beta ,\ t\in
S^1,$$
reduces to 
$$(f'(\tau )\partial _t+\mu (\tau )(\alpha -\beta ))u(t,\tau )=v(t,\tau
),$$ and has a unique solution $u$ for any given (smooth) $v$ when $\alpha
\ne
\beta $. If $\alpha =\beta $, we have a solution which is unique up to a
$\tau $-dependent constant, provided that we replace $v(t,\tau )$ by
$v(t,\tau )-\langle v(\cdot ,\tau )\rangle $.

\par It follows in the orientable case, that if $v(t,\tau ,x,\xi )$ is a
\hm{} \pol{} of degree $j$ in $x,\xi $ with analytic \coeff{}s depending
on $(t,\tau )\in S^1\times ]-\epsilon _0,\epsilon _0[$, then we can find
$u(t,\tau ,x,\xi )$ of the same type such that 
\ekv{1.20}
{(f'(\tau )\partial _t+\mu (\tau )(x\partial _x-\xi \partial _\xi
))u(t,\tau ,x,\xi )=v(t,\tau ,x,\xi )-[v](\tau ,x,\xi ),}
where 
$$[v](\tau ,x,\xi )=\sum_{\alpha}\langle v_{\alpha ,\alpha
}(\cdot ,\tau )\rangle (x\xi )^\alpha ,\ v=\sum_{\vert \alpha
+\beta\vert =j } v_{\alpha ,\beta }(t,\tau )x^\alpha \xi ^\beta .$$

\par In the non-orientable case, $x,\xi $ are anti-\per{} in $t$, so the
general form of a \fu{} of $(t,\tau )\in S^1\times ]-\epsilon _0,\epsilon
_0[$ with values in the space of $j$-\hm{} \pol{}s in $(x,\xi )$ is 
\ekv{1.21}
{v(t,\tau ,x,\xi )=\sum_{\vert \alpha +\beta \vert =j}e^{i(\alpha -\beta
)t/2}v_{\alpha ,\beta }(t,\tau ) x^\alpha \xi ^\beta .}
Using that $\mu (\tau )>0$ is real, we get the same result as in the
orientable case for the solvability of (1.20).

\par In both cases, we combine the solvability of (1.20) with (1.19) and
see that there is a sequence of $G_3,G_4,...$ as above, so that at the
level of formal Taylor series in $(x,\xi )$:
\ekv{1.22}
{p\circ \exp H_{G_3}\circ \exp H_{G_4}\circ ... \, =f(\tau )+\mu (\tau
)x\xi +q(\tau ,x\xi ),}
where 
\ekv{1.23}
{q(\tau ,x\xi )=\sum_{\alpha =2}^\infty  q_\alpha (\tau )(x\xi )^\alpha }
is resonant in the sense that it is a \fu{} of $\tau $ and $x\xi $ only. 

\par Since we do not wish to consider convergence questions here (even
though the Birkhoff series are likely to converge in the present
2-dimensional case), we stop at some high but finite order and write
\ekv{1.24}
{p\circ \kappa _N=f(\tau )+\mu (\tau )x\xi +q^{(N)}(\tau ,x\xi )+{\cal
O}((x,\xi )^{N+1}),}
where 
\ekv{1.25}
{
q^{(N)}(\tau ,x,\xi )=\sum_{4\le 2\alpha \le N}q_\alpha (\tau )(x\xi
)^\alpha ,\ \kappa _N=\exp H_{G_3}\circ ...\circ\exp H_{G_N}. }
Later, we shall review how to reduce the original \op{} (1.1) to an
$h$-\pop{} $P$ with the new principal symbol 
\ekv{1.26}
{p=f(\tau )+\mu (\tau )x\xi +q^{(N)}(\tau ,x\xi )+{\cal O}((x,\xi
)^{N+1}). }
We recall how the lower order symbols in $P$ can be simplified by
conjugation with elliptic \pop{}s: 

\par Write 
\ekv{1.27}
{
P(t,\tau ,x,\xi ;h)=p(t,\tau ,x,\xi )+hp_1(t,\tau ,x,\xi )+...,
}
where in general, we identify symbols with their $h$-Weyl quantizations.
If $A(t,hD_t,x,hD_x;h)$ is an $h$-\pop{} of order $-m\le 0$ in $h$ and with
symbol $A=h^ma_m+h^{m+1}a_{m+1}+... $, with $a_j(t,\tau ,x,\xi )$ smooth
in some fixed domain, we use that on the operator level
\ekv{1.28}
{e^{-iA}Pe^{iA}=e^{-i{\rm ad}_A}P=\sum_{k=0}^\infty {1\over k!}(-i{\rm
ad}_A)^kP,}
to see that $e^{-iA}Pe^{iA}-P$ is of order $-(m+1)$ and has the leading
symbol
\ekv{1.29}
{h^{m+1}\{ p,a_m\} .}
If we first let $m=0$, we can repeat the discussion above for the $G_j$,
and see that we can choose $a_0=a_0^{(N)}$, so that 
\ekv{1.30}
{
\{ p,a_0^{(N)}\} -p_1=b^{(N)}(\tau ,x\xi )+{\cal O}((x,\xi )^{N+1}).
}
Proceeding similarly with the lower order symbols, we see that we can
find $A=A^{(N)}$ with symbol
$a_0^{(N)}+ha_1^{(N)}+..+h^{N-1}a_{N-1}^{(N)}$, such that 
\ekv{1.31}
{
e^{-iA}Pe^{iA}=P^{(N)}+R_{N+1}
}
where $P^{(N)}$ has the symbol $P^{(N)}(\tau ,x\xi ;h)$ and
$R_{N+1}(t,\tau ,x,\xi )={\cal O}((h,x,\xi )^{N+1})$. Moreover,
\ekv{1.32}
{P^{(N+1)}(\tau ,x\xi;h) -P^{(N)}(\tau ,x\xi ;h)={\cal O}((h,x,\xi
)^{N+1}).} 

\par Summing up the discussion so far, we have
\medskip
\par\noindent \bf \Prop{} 1.1. \it
Make the assumptions above, in particular that $P$ has analytic \coef{}s
and satisfies (1.7), (1.8).

\par There exists an analytic \ctf{} $\kappa :{\rm neigh\,}(\{ \tau
=x=\xi =0\} ,T^*S^1\times {\bf R})\to {\rm neigh\,}(\gamma (E_0),T^*{\bf
R}^2)$, single-valued in the orientable case and otherwise double-valued
with $\kappa (t-2\pi ,x,\xi )=\kappa (t,-x,-\xi )$, such that
\ekv{1.33}
{p\circ \kappa =f(\tau )+\mu (\tau )x\xi +q^{(N)}(\tau ,x\xi )+{\cal
O}((x,\xi )^{N+1}),}
with $q^{(N)}$
as in (1.25), and $\tau =x=\xi =0$ corresponding to $\gamma (E_0)$. Here
$N\ge 1$ is any fixed integer.

\par Let $S=\int_{\gamma (E_0)}\xi dx$ be the action of $\gamma (E_0)$. Let
$L^2_S(S^1\times {\bf R})$ be the space of locally square integrable
\fu{}s $u(t,x)$ on
${\bf R}\times {\bf R}$ with 
$$\Vert u\Vert ^2:=\int \int_0^{2\pi }\vert u(t,x)\vert ^2dtdx 
<\infty ,$$
$u(t-2\pi ,x)=e^{iS/h}u(t,x)$ in the orientable case, $u(t-2\pi
,x)=e^{iS/h}u(t,-x)$ in the non-orientable case. Then there exists an
analytic unitary \fop{} $U:L^2_S(S^1\times {\bf R})\to L^2({\bf R}^2)$
associated to $\kappa $, and microlocally defined near $\tau =x=\xi =0$ in
$T^*(S^1\times {\bf R})$ such that (microlocally):
\ekv{1.34}
{PU=U(P^{(N)}+R_{N+1}),}
where $P^{(N)},R_{N+1}$ have the Weyl symbols $P^{(N)}(\tau ,x\xi ;h)\sim
\sum p_j^{(N)}(\tau ,x\xi )h^j$, $R_{N+1}(t,\tau ,x,\xi ;h)={\cal
O}((h,x,\xi )^{N+1})$, satisfying (1.32).
\rm
\medskip

\par Here (as in [MeSj]), we define our \fop{} $U$ on the FBI-Bargman
\tf{} side.

\bigskip

\centerline{\bf 2. Resonances.}
\medskip
\par It follows from [GeSj] (see also [Sj2]) that we can find  a smooth
function $G(x,\xi )\in C^\infty ({\bf R}^4;{\bf R})$ which is an escape
\fu{} in the sense of [HeSj] and satisfies:
\ekv{2.1}
{
H_pG>0\hbox{ in }p^{-1}(]E_0-\epsilon _0,E_0+\epsilon _0[)\setminus
\Gamma , }
\ekv{2.2}
{H_pG\sim {\rm dist\,}((x,\xi ),\Gamma )^2,\hbox{ near }\Gamma
\hbox{ in }p^{-1}(]E_0-\epsilon _0,E_0+\epsilon _0[).}
In the coordinates $(t,\tau ,x,\xi )$, where $p$ is reduced to the \rhs{}
of (1.24), we may even assume that 
\ekv{2.3}
{G(t,\tau ,x,\xi )={1\over 2}(x^2-\xi ^2).}
Recall from the [HeSj]-theory that we have an associated IR-\mfld{}
$\Lambda _\delta \subset {\bf C}^4$, defined by 
$$\Lambda _\delta =\{ (x,\xi )=\exp (i\delta H_G)(y,\eta );\, (y,\eta
)\in {\bf R}^4\} ,$$
for $0<\delta \ll 1$. (Strictly speaking, the above representation of
$\Lambda _\delta $ is valid only where $G$ is analytic and with $G$
denoting also the \hol{} extension. Elsewhere, we let $G$ denote a
suitable almost \hol{} extension and take $(x,\xi )=\exp (\delta H^{\Im
\sigma }_{\Re G})(y,\eta )$, where $H_{\Re G}^{\Im \sigma }$ denotes the
Hamilton field of $\Re G$ \wrt{} the real symplectic form $\Im \sigma $.)
Using (2.3) we see that in the special coordinates used there,
$\Lambda _\delta
$ is given by:
\ekv{2.4}
{t,\tau \in{\bf R},\ \cases{x=(\cos \delta )y-i(\sin \delta )\eta \cr \xi
=-i(\sin \delta )y+(\cos \delta ) \eta },\,\, y,\eta \in{\bf R}.}
This can be written,
\ekv{2.5}
{
t,\tau \in{\bf R},\ \xi ={2\over i}{\partial \Phi _\delta (x)\over
\partial x},\ \Phi _\delta (x)={1\over 2}(\cot \delta )(\Im x)^2+{1\over
2}(\tan \delta )(\Re x)^2. }

\par Notice that $\Phi _{\pi /4}(x)={1\over 2}\vert x\vert ^2$ and that the
corresponding IR-\mfld{} is given by $t,\tau \in {\bf R}$, $\xi
=-i\overline{x}$. Following an argument from [KaKe], we look for a new
IR-\mfld{} $\Lambda $ which coincides with $\Lambda _{\pi /4}=\Lambda
_{\Phi _{\pi /4}}$ near
$\gamma _{E_0}$ and with $\Lambda _\delta $ outside a small \neigh{} of $\gamma
_{E_0}$. First we notice that if $q_1(x),q_2(x)$ are \st{} convex
quadratic forms on ${\bf R}^2$, then we can find a smooth \st{} convex
function $\phi (x)$, with $\phi (x)=q_1(x)$ near $0$ and with $\phi
(x)=q_2(x)$ outside $V$, where $V$ is any given \neigh{} of $0$. Apply
this with $q_1(x)=\Phi _{\pi /4}(x)$, $q_2(x)=\Phi _\delta (x)$, then
replace
$\phi (x)$ by $\alpha ^2\phi (x/\alpha )$, $0<\alpha \ll 1$, in order
to decrease the \neigh{} of $0$, where $\phi \ne \Phi _\delta $ even
further while keeping $\phi $ \bdd{} in $C^2 $. Using the strict
convexity of $\phi (x)$ and the fact that $\phi (0)=0$ is  global
minimum, we see that 
\ekv{2.6}
{
-\Im (x{2\over i}{\partial \phi \over \partial x})=\Re (2x{\partial \phi
\over \partial x})=\langle x,\nabla \phi \rangle _{{\bf R}^2}\sim \vert
x\vert ^2, }
\ufly{} in $\alpha $ but not in $\delta $. Combining this with (1.15), we
see that for $x\in{\bf C}$, $\vert x\vert <1/{\cal O}(1)$ \indep{}ly of
$\alpha ,\delta $:
\ekv{2.7}
{
\Im p(t,\tau ,x,{2\over i}{\partial \phi \over \partial x})\sim -x^2,
}
\ekv{2.8}
{
\vert \Re p(t,\tau ,x,{2\over i}{\partial \phi \over \partial
x})-E_0\vert \ge {1\over C}\vert \tau -\tau _0\vert -C\vert x\vert ^2,\
t,\tau \hbox{ real,} }
where $\tau _0$ is the value with $f(\tau _0)=E_0$. We may assume without
loss of generality that $\tau _0=0$.

\par As a first attempt to modification of $\Lambda _\delta $, we
consider 
$\{ (t,\tau ,x,\xi );\, t,\tau \in{\bf R},\, \xi ={2\over i}{\partial
\phi \over \partial x}\}$, but we need the modification to coincide with
$\Lambda _\delta $, not only for $\vert (x,\xi )\vert $ outside a
\neigh{} of $0$, but also for $\vert \tau \vert $ away from $0$. In order
to do so, we change the representation of $\Lambda _\delta $ in
accordance with the application of the Bargman \tf{},
\ekv{2.9}
{
u(t)\mapsto h^{-3/4}\int e^{-(t-s)^2/2h}u(s)ds,
}
whose corresponding \canform{} sends $S_t^1\times {\bf R}_\tau $ into $\{
(t,\tau )\in (S^1+i{\bf R})\times {\bf C};\, \tau ={2\over i}{\partial
\over \partial t}({1\over 2}(\Im t)^2)\, (=-\Im t)\} $. Since (2.9)
is a convolution \op{}, $f(hD_t)$ keeps the same shape after the
transformation, and so does the principal symbol $p$ in (1.15). 

\par After applying the \canform{} associated to the \op{} (2.9),
$\Lambda _\delta $ becomes
\ekv{2.10}
{
(\tau ,\xi )={2\over i}{\partial \over \partial (t,x)}\widetilde{\Phi
}_\delta (t,x),\ \widetilde{\Phi }_\delta (t,x)={1\over 2}(\Im t)^2+\Phi
_\delta (x), }
and the symbol $p$ is still of the form (1.24). Let $\delta >0$ be small
but fixed, let $\chi (\Im t)$ be a standard cutoff around $\Im t=0$ and
consider the "intermediate" weight
\ekv{2.11}
{\widetilde{\Phi }(t,x)={1\over 2}(\Im t)^2+\chi (\Im t)\phi (x)+(1-\chi
(\Im t))\Phi _\delta (x),}
which satisfies
\ekv{2.12}
{
\widetilde{\Phi }(t,x)=\widetilde{\Phi }_\delta (t,x),\hbox{ for }\vert
\Im t\vert +\vert x\vert \ge {\rm Const.}, }
\ekv{2.13}{\widetilde{\Phi }(t,x)=\widetilde{\Phi }_{\pi /4}(t,x),\hbox{
for }\vert
\Im t\vert +\vert x\vert \le {\rm const.}.}
Notice that the mixed derivative $\nabla _t\nabla _x \widetilde{\Phi
}(t,x)$ has its support in $\vert x\vert ={\cal O}(\alpha )$ and is
therefore ${\cal O}(\alpha )$. Consequently $\widetilde{\Phi }$ is \st{}
convex as a \fu{} of $\Im t,\, x$ and in particular \stpsh{}. 

\par We have 
$$\eqalign{
{2\over i}{\partial \widetilde{\Phi }\over \partial x}&=\chi (\Im
t){2\over i}{\partial \phi \over \partial x}+(1-\chi (\Im t)){2\over
i}{\partial \Phi _\delta (x)\over \partial x},\cr 
{2\over i}{\partial \widetilde{\Phi}\over \partial t}&=-\Im t-(\phi
(x)-\Phi _\delta (x))\chi '(\Im t). }$$
Notice that the last expression is real, so in view of the properties of
$\phi $, $\Phi _\delta $, and the fact that 
\ekv{2.14}
{
p=f(\tau )+\mu (\tau )x\xi +{\cal O}((x,\xi )^4),
}
we get with $\Lambda =\Lambda _{\widetilde{\Phi }}$: $\xi ={2\over
i}{\partial
\widetilde{\Phi }\over \partial x},\, \tau ={2\over i}{\partial
\widetilde{\Phi }\over \partial t}$:
\ekv{2.15}
{
{\Im p_\vert}_{\Lambda _{\widetilde{\Phi }}}=\Im p(t,{2\over i}{\partial
\widetilde{\Phi }\over \partial t},x,{2\over i}{\partial \widetilde{\Phi
}\over \partial x})\sim -\vert x\vert ^2 ,}
\ekv{2.16}
{{\Re p_\vert}_{\Lambda _{\widetilde{\Phi }}}=f(-\Im t)+{\cal O}(\vert
x\vert ^2).}
In particular, if we assume that $\gamma (E_0)$ corresponds to $\tau =0$,
so that $f(0)=E_0$, then
$$\vert {(p-E_0)_\vert}_{\Lambda _{\widetilde{\Phi }}}\vert \sim \vert \Im
t\vert +\vert x\vert ^2,$$
so the study of the spectrum of $p$ near $0$ depends essentially only on
the behaviour of our \op{} in a \neigh{} of $\tau =0$, $x=\xi =0$, where
$\Lambda _{\widetilde{\Phi }}=\Lambda _{\Phi _{\pi /4}}$.

\par In such a \neigh{} we identify $\Lambda _{\Phi _{\pi /4}}$ linearly
and symplectically with
$T^*(S^1_t\times {\bf R}_y)$ (actually $T^*M$ in the non-orientable case,
where $M$ is defined after \Prop{} 2.1) in such a way that
${p_\vert}_{\Lambda _{\Phi _{\pi /4}}}$ becomes
$$\eqalignno{p&=f(\tau )+\mu (\tau ){1\over 2i}(y^2+\eta ^2)+q^{(N)}(\tau
,{1\over 2i}(y^2+\eta ^2))+{\cal O}((y,\eta )^{N+1}),&(2.17)\cr
q^{(N)}&=\sum_{4\le 2\alpha \le N}q_{\alpha (\tau )}({1\over 2i}(y^2+\eta
^2))^\alpha .} $$
More explicitly, on $\Lambda _{\Phi _{\pi /4}}$ we have $\xi
=-i\overline{x}$, $d\xi \wedge dx=-id\overline{x}\wedge dx$ and we can
take $\eta =\sqrt{2}\, \Re x$, $y=\sqrt{2}\,\Im x$. From now on, we write
$x,\xi $ instead of $y,\eta $ for the new coordintates. This identification
can be quantized by means of a metaplectic FBI
\tf{}, so applying (1.31) (rather than redoing the arguments leading to
that equation) we may assume after conjugating with an elliptic $h$-\pop{},
that we have 
\ekv{2.18}
{P(t,\tau ,x,\xi ;h)=P^{(N)}(\tau ,{1\over 2i}(x^2+\xi ^2);h)+R_{N+1},}
where $R_{N+1}={\cal O}((h,x,\xi )^{N+1})$ and changing the lower order
terms here, that on the \op{} level:
\ekv{2.19}
{
P^w(t,hD_t,x,hD_x;h)=P^{(N)}(hD_t,{1\over
2i}(x^2+(hD_x)^2);h)+R_{N+1}(t,hD_t,x,hD_x;h). }
Here the obvious analogues of (1.32) hold and the $h$-principal symbol
of $P$ is equal to $p$ in (2.17).

\par Summing up the discussion so far, we have
\medskip
\par\noindent \bf \Prop{} 2.1. \it
Let $\delta >0$ be small enough. We can find a closed IR-\mfld{} 
$\Lambda $ which coincides with $\Lambda _\delta $ outside an \ably{}
small \neigh{} of $\gamma (E_0)$ and contains $\gamma (E)$ for $E\in{\rm
neigh\,}(E_0,{\bf R})$, such that 
\smallskip
\par\noindent 1) ${(p-E)_\vert}_{\Lambda }\ne 0$ on $\Lambda \setminus
\gamma (E)$, $E\in{\rm neigh\,}(E_0,{\bf R})$.\smallskip
\par\noindent 2) In the coordinates $(t,\tau ,x,\xi )$ of \Prop{} 1.1,
and after applying the \ctf{} associated to (2.9), so that $T^*S^1$
becomes $\{ (t,\tau );t\in S^1+i{\bf R},\, \tau \in {\bf C},\, \tau =-\Im
t\}$, $\Lambda _\delta $ takes the form (2.10) and $\Lambda $ becomes
$\Lambda _{\widetilde{\Phi }}$, where $\widetilde{\Phi }(t,x)$ is
\stpsh{}, satisfying (2.12,13) and 
$$\vert {(p-E)_\vert }_{\Lambda _{\widetilde{\Phi }}}\vert \sim \vert \Im
t-f^{-1}(E)\vert +\vert x\vert ^2,\ E\in{\rm neigh\,}(E_0,{\bf R}).$$
\smallskip
\par\noindent 3) Let $T$ be the \op{} (2.9), and let $H^{{\rm
loc}}_{\widetilde{\Phi },S}(\Omega )$, be the space of \hol{} \fu{}s
$u(t,x)$ on $\Omega $ equipped with the seminorms 
$$\Vert u\Vert _{H_{\widetilde{\Phi },S}(K)}^2=\int _K \vert u(t,x)\vert ^2
e^{-2\widetilde{\Phi }(t,x)/h}L(dtdx)$$ with the Floquet condition 
$u(t-2\pi ,x)=e^{iS/h}u(t,\pm x)$, with $+$ in the orentiable case and
$-$ otherwise. Then with $U$
as in \Prop{} 1.1, $\widetilde{P}:=(T\otimes 1)\circ U^{-1})\circ P\circ
U\circ (T^{-1}\otimes 1)$ is a well-defined $h$-\pop{} acting on
$H_{\widetilde{\Phi },S}^{\rm loc}(\Omega )$, of the form
\ekv{2.20}
{\widetilde{P}=P^{(N)}(hD_t,{1\over 2}(xhD_x+hD_xx);h)+R_{N+1}, }
with $P^{(N)}$, $R_{N+1}$ as above.\smallskip
\par\noindent 4) In the region where $\widetilde{\Phi }=\Phi _{\pi /4}$,
we can identify $\widetilde{P}$ with $P^{(N)}(hD_t,{1\over
2i}(x^2+(hD_x)^2);h)+\widetilde{R}_{N+1}$ acting on $L_{S,{\rm
loc}}^2(S^1\times {\bf R})$. 
\rm\medskip
\medskip

\par For future use, we introduce some geometrical objects that we have
already (implicitly) encountered, and make a very short review of how they
fit to the definition of \res{}s in the [HeSj]-theory. Let $M=S^1\times
{\bf R}$ in the orientable case and in the non-orientable case, let $M={\bf
R}_t\times {\bf R}_x/\sim$, where $(t,x)\sim (t-2\pi ,-x)$. In both
cases, we use $(t,\tau ,x,\xi )$ as local coordinates on $T^*M$. Then in
\Prop{} 1.1, we have $\kappa :\,{\rm neigh\,}(\tau =x=\xi =0,T^*M )\to
{\rm neigh\,}(\gamma (E_0),T^*{\bf R}^2)$. The Weyl symbol of $U^{-1}PU$ in
(1.34) is defined on $T^*M$. A natural complexification $\widetilde{M}$ of
$M$ is given by $(S^1+i{\bf R})\times {\bf C}$ in the orientable case,
and by ${\bf C}_t\times {\bf C}_x/\sim$, in the non-orientable case, where
$\sim$
is defined as above, now for complex $t,x$. If $\Omega \subset {\bf
C}_t\times {\bf C}_x$ is open and invariant under $(t,x)\mapsto (t,-x)$
and $(t,x)\mapsto (t-2\pi ,x)$, then we can view $\Omega $ as a subset of
$\widetilde{M}$, and $H^{{\rm loc}}_{\widetilde{\Phi },S}(\Omega )$ is
defined in \Prop{} 2.1.

\par As in [KaKe], we can define the Hilbert space $H(\Lambda )$ in the
spirit of [HeSj], to be the space of $u\in H(\Lambda _\delta )$ equipped
with the norm 
$$\Vert (1-\chi _2)T_{{\rm HS}}u\Vert _{L^2(T^*{\bf R}^2,e^{-2H_\delta
/h}dxd\xi )}+\Vert \chi _1(t,x)(T\otimes
1)U^{-1}u\Vert_{L^2(e^{-2\widetilde{\Phi }}L(dxdt))} ,$$ where $\chi _1\in
C_0^\infty (\widetilde{M})$ is equal to 1 in a \neigh{} of $\Im t=x=0$,
containing the region where
$\widetilde{\Phi }$ differs from $\Phi _\delta $, and $\chi _2\in
C_0(T^*{\bf R}^2)$ is equal to 1 near $\gamma (E_0)$ and is essentially
the same function as $\chi _1$, after suitable identifications of domains.
$H_\delta $ is the weight appearing in the definition ([HeSj]) of
$H(\Lambda _\delta )$ (for which the corresponding norm would have been
given by only the first term in the expression above, without the factor
$1-\chi _2$), and $T_{{\rm HS}}$ is the correponding global FBI-\tf{}. The
resonances near $E_0$ are then the \ev{}s of $P$, viewed as an \op{} in
$H(\Lambda )$, (with domain $H(\Lambda ,m)$, where $m\ge 1$ is the order
\fu{} associated to $P$, in the the most standard case, $m(x,\xi
)=\langle \xi \rangle ^2$).

\par If $u_{\ell}(x)$, $\ell =0,1,2,..$ denote the normalized \ef{}s of
${1\over 2}(x^2+(hD)^2)$ with \ev{}s $(\ell +{1\over 2})h$, we have the ON
system of \ef{}s to $P^{(N)}(hD_t,{1\over 2i}(x^2+(hD_x)^2);h)$:
\ekv{2.21}
{v_{k,\ell}={1\over \sqrt{2\pi }}e^{{i\over h}(kh-S/2\pi
)t}u_{\ell}(x),\ k\in{\bf Z}, \ell\in {\bf N},}
in the orientable case, and 
\ekv{2.22}
{v_{k,\ell}={1\over \sqrt{2\pi }}e^{{i\over h}((k+{\ell\over 2})h-S/2\pi
)t}u_{\ell}(x),\ k\in{\bf Z}, \ell\in {\bf N},}
in the non-orientable case
with the corresponding \ev{}s:
\eekv{2.23}
{
P^{(N)}(hk-{S\over 2\pi },{1\over i}(\ell +{1\over 2})h;h),
} 
{P^{(N)}(h(k+{\ell\over 2})-{S\over 2\pi },{1\over i}(\ell +{1\over
2})h;h),}
in the orientable and non-orientable cases respectively. Notice that the
\ef{}s (2.21,22) are
\ml{}ly concentrated to the region in phase space, where $\tau \approx
hk-S/2\pi
$, ${1\over 2}(x^2+\xi ^2)\approx (\ell +{1\over 2})h$. (Such remarks were
used in [Sj4].) Combining such arguments with the reduction above and some
arguments of [MeSj] we get
\medskip
\par\noindent \bf \Prop{} 2.2. \it Fix $\delta \in ]0,1[$ and $N\in{\bf
N}$ with $h^{(N+1)\delta /2}\ll h$, i.e. with $(N+1)\delta /2>1$. Then
the \res{}s of $P$ in (1.1) in the rectangle $]-E_0-\epsilon
_0,E_0+\epsilon _0[-i[0,h^\delta [$ are simple and given by 
\eekv{2.24}{P^{(N)}(hk-S/2\pi ,{1\over i}(\ell +{1\over 2})h;h)+{\cal
O}(((\ell +{1\over 2})h)^{{(N+1)\over 2}}),}
{P^{(N)}(h(k+{\ell\over 2})-S/2\pi ,{1\over i}(\ell +{1\over 2})h;h)+{\cal
O}(((\ell +{1\over 2})h)^{{(N+1)\over 2}}),}
respectively in the orientable and the non-orientable cases,
for $k\in {\bf Z}$, $\ell \in{\bf N}$, $\ell h={\cal O}(h^\delta
)$.\rm\medskip

\par Now use the comptatibility property (1.32) and let $\widetilde{P}$
satisfy 
\ekv{2.25}
{\widetilde{P}(\tau ,\iota ;h)=P^{(N)}(\tau ,\iota ;h)+{\cal O}((h,\iota
^{1/2})^{N+1})}
fo all $N$, 
\ekv{2.26}
{\widetilde{P}(\tau ,\iota ;h)\sim \widetilde{p}(\tau ,\iota
)+h\widetilde{p}_1(\tau ,\iota )+...\,.}
Then \Prop{} 2.2 shows that the \res{}s in the rectangle $]E_0-\epsilon
_0,E_0+\epsilon _0[-i[0,h^\delta [$ are simple and of the form,
\eekv{2.27}
{\widetilde{P}(hk-{S\over 2\pi }, {1\over i}(\ell +{1\over 2})h;h)+{\cal
O}(h^\infty ),\ k\in{\bf Z},\ell\in{\bf N}\hbox{ (orientable case)}, }
{\widetilde{P}(h(k+{\ell\over 2})-{S\over 2\pi }, {1\over i}(\ell +{1\over
2})h;h)+{\cal O}(h^\infty ),\ k\in{\bf Z},\ell\in{\bf N}\hbox{
(non-orientable case).} }

\par Next we study the \ev{}s $z$ with $h^\delta <-\Im z\le \epsilon _1$
for some small $\epsilon _1>0$ \indep{} of $h$. If $-\Im z\sim \epsilon
$, the \ef{}s will be localized to a region with $\vert (x,\xi )\vert
\sim \epsilon ^{1/2}$, so we make the change of variables 
$$x=\epsilon ^{1/2}\widetilde{x},\,
hD_x=\epsilon ^{1/2}\widetilde{h}D_{\widetilde{x}},\, \widetilde{h}={h\over
\epsilon }.$$ Then we also have $hD_t=\epsilon \widetilde{h}D_t$, and
$P^w$ in (2.19) becomes
\ekv{2.28}
{P=P^{(N)}(\epsilon \widetilde{h}D_t,\epsilon {1\over
2i}(\widetilde{x}^2+(\widetilde{h}D_{\widetilde{x}})^2);h)+R_{N+1}(t,\epsilon
\widetilde{h}D_t,\epsilon ^{1\over 2}\widetilde{x}, \epsilon ^{1\over
2}\widetilde{h}D_{\widetilde{x}};\epsilon \widetilde{h}),}
to be studied in the region where where $\vert
(\widetilde{x},\widetilde{\xi })\vert \sim 1$. In this region, the symbol
of the remainder term in (2.28) is ${\cal O}((\epsilon
\widetilde{h})^{N+1}+\epsilon ^{N+1\over 2})$. 
(Recall that (2.19) describes the original \op{} acting in $H(\Lambda )$,
viewed in a \neigh{} of $\gamma (E_0)$. Correspondingly (2.19) is \ml{}ly
defined near $\tau =x=\xi $ in $T^*M$ and acts on $L^2_S(M)$.) 

\par Let 
$$P^{(N)}(\tau ,\sigma ;h)\sim\sum_{j=0}^\infty  p_j^{(N)}(\tau ,\sigma
)h^j,$$
so that 
$$
P^{(N)}(\epsilon \tau ,\epsilon \sigma ;h)\sim \sum_{j=0}^\infty
p_j^{(N)}(\epsilon \tau ,\epsilon \sigma )\epsilon
^j\widetilde{h}^j=\epsilon ({1\over \epsilon }p_0^{(N)}(\epsilon \tau
,\epsilon \sigma )+\sum_{j=1}^\infty p_j^{(N)}(\epsilon \tau
,\epsilon \sigma )\epsilon ^{j-1}\widetilde{h}^j ).
$$
Here 
$$
{1\over \epsilon }(p_0^{(N)}(\epsilon \tau ,\epsilon \sigma )-
E_0)={1\over
\epsilon }(f(\epsilon \tau )-E_0)+{1\over i}\mu (\epsilon \tau )\sigma
+{\cal O}(\epsilon ).
$$

\par Notice that 
$$\Lambda _{E,F}^0:=\{(t,\tau ,x,\xi )\in T^*M;\, {1\over \epsilon
}(f(\epsilon \tau )-E_0)=E,\, \mu (\epsilon \tau )\sigma =-F\}$$
is a family of Lagrangian torii for $E,F-1\in {\rm neigh\,}(0,{\bf R})$
which form an analytic foliation of phase-space and depend analytically
on $\epsilon \in{\rm neigh\,}(0,{\bf R})$. This means that the
$\widetilde{h}$-\pop{} (2.28) fullfills the assumptions of [MeSj] (with the
slight difference that we are now on a 2-dimensional analytic \mfld{}
$M$, rather than ${\bf R}^2$). The geometric main result of [MeSj] tells
us that there is a corresponding foliation of $T^*\widetilde{M}$ into
complex Lagrangian torii, $\Lambda _{E,F}$, for $E,F\in{\rm
neigh\,}(0,{\bf C})$, which is $\epsilon $-close to the complexification
of the family $\Lambda ^0_{E,F}$, and such that ${1\over \epsilon
}(p_0-E_0)$ is constant $=E+iF+{\cal O}(\epsilon )$ on each $\Lambda
_{E,F}$, where $p_0$ denotes the principal symbol of the \op{} (2.28). 

\par The corresponding spectral result of [MeSj] is applicable also
(with some simple and straight forward modifications in the definition
of a certain global Grushin \pb{} for the original operator $P$ acting
on $H(\Lambda )$ in the proof). We conclude that (the original \op{})
$P$ has the
\res{}s
$g(k\widetilde{h}-{S\over 2\pi
\epsilon },(\ell +{1\over 2})\widetilde{h},\epsilon ;\widetilde{h})$ in
the orientable case and  
$g((k+{\ell\over 2})\widetilde{h}-{S\over 2\pi
\epsilon },(\ell +{1\over 2})\widetilde{h},\epsilon ;\widetilde{h})$ in
the non-orientable case, where 
\ekv{2.29}
{
g(\tau ,\sigma ,\epsilon ;\widetilde{h})\sim \sum_{j=0}^\infty
g_j(\tau ,\sigma ,\epsilon )\widetilde{h}^j. }
(Notice here that we can \tf{} ${1\over
2i}(\widetilde{x}^2+(\widetilde{h}D_{\widetilde{x}})^2)$, into ${1\over
2}(x\widetilde{h}D_x+\widetilde{h}D_xx)$ by standard Bargman \tf{}, and
that the last \op{} becomes ${1\over i}h(D_s+{1\over 2})$, if we make the
change of variables $x=e^{is}$.)

\par The corresponding \ev{}s of $P^{(N)}(\epsilon
\widetilde{h}D_t,{\epsilon \over 2i}\widetilde{h}(D_s+{1\over 2});h)$ are 
$g^{(N)}(k\widetilde{h}-{S\over 2\pi \epsilon },(\ell +{1\over
2})\widetilde{h},\epsilon ;\widetilde{h})+{\cal O}(h^\infty )$ (in the
orientable case and with the usual modification in the non-orientable
case), with 
\ekv{2.30}
{
g^{(N)}(\tau ,\sigma ,\epsilon ;\widetilde{h})\sim
\sum_{j=0}^\infty  g_j^{(N)}(\tau ,\sigma ,\epsilon )\widetilde{h}^j, }
given by $g^{(N)}(\tau ,\sigma ,\epsilon ;\widetilde{h})=P^{(N)}(\epsilon
\tau ,\epsilon \sigma ;\epsilon \widetilde{h})$. Because of the estimate on
$R_{N+1}$, we know from the proof in [MeSj] that 
\ekv{2.31}
{
g_j(\tau ,\sigma ,\epsilon )-g_j^{(N)}(\tau ,\sigma ,\epsilon )={\cal
O}(\epsilon ^{N+1\over 2}),\hbox{ for }j\le N. }

\par Writing
$$P^{(N)}(\tau ,\sigma ;h)=\sum p_j^{(N)}(\tau ,\sigma )h^j,$$
we see that 
\ekv{2.32}
{
g_j^{(N)}(\tau ,\sigma ,\epsilon )=p_j^{(N)}(\epsilon \tau ,\epsilon \sigma
)\epsilon ^j. }

\par The choice of $\epsilon $ is not unique. Substituting $(\epsilon
,\widetilde{h})\mapsto (\mu \epsilon ,\widetilde{h}/\mu )$, with $\mu
\sim 1$, will not affect  the \ev{}s of (2.28) (which is only a rewriting
of the $\epsilon $-\indep{} \op{} (2.19)), so
$$g(k{\widetilde{h}\over \mu }-{S\over 2\pi \mu
\epsilon },(\ell +{1\over 2}){\widetilde{h}\over \mu },\mu \epsilon
;{\widetilde{h}\over \mu })=g(k\widetilde{h}-{S\over 2\pi
\epsilon },(\ell +{1\over 2})\widetilde{h},\epsilon ;\widetilde{h}),$$
(with the usual modification in the non-orientable case)
and as in section 7 of [MeSj], we conclude that
\ekv{2.33}
{
g_j({\tau \over \mu },{\sigma \over \mu },\mu \epsilon )=\mu ^jg_j(\tau
,\sigma ,\epsilon ). }
The same relation holds for $g_j^{(N)}$, as can also be seen directly
from (2.32).

\par Now recall that $g_j(\tau ,\sigma ,\epsilon )$ is defined for $\tau
={\cal O}(1/\epsilon )$,
$\sigma \sim 1$. Using (2.33), we define $g_j(\tau ,\sigma ,1)$ by 
\ekv{2.34}
{g_j(\tau ,\sigma ,1)= \mu ^{-j}g_j({\tau \over \mu },{\sigma \over \mu
},\mu ),}
when $0<\sigma \ll 1$, $\tau ={\cal O}(1)$, by taking $\mu \sim \sigma $.
Similarly, we have $g_j^{(N)}(\tau ,\sigma ,1)=p_j^{(N)}(\tau ,\sigma )$,
which is analytic in $\tau ={\cal O}(1)$, $0\le \sigma \ll 1$. (2.31)
implies that 
\ekv{2.35}
{g_j(\tau ,\sigma ,1)=p_j^{(N)}(\tau ,\sigma )+{\cal O}(\sigma
^{{N+1\over 2}-j}),\ j\le N.}
Using that we have estimates of the type (2.31) also for the
derivatives, we conclude that $g_j(\tau ,\sigma ,1)$ is smooth down to
$\sigma =0$. 

\par The formula for the \res{}s of $P$ prior to (2.29) together with
(2.29) shows that we have the \res{}s:
\ekv{2.36}
{
\sum_{j=0}^\infty  g_j({kh\over \epsilon }-{S\over 2\pi \epsilon },(\ell
+{1\over 2}){h\over \epsilon },\epsilon )\epsilon
^{-j}h^j=\sum_{j=0}^\infty g_j(kh-{S\over 2\pi},(\ell +{1\over
2})h,1)h^j,} where $\epsilon \ge h^\delta $, $\ell , k$ are chosen so that
$(\ell+{1\over 2})h/\epsilon \sim 1$ and $kh-S/2\pi ={\cal O}(1)$.

\par We get, using also \Prop{} 2.2:
\medskip
\par\noindent \bf\Th{} 2.3. \it Same assumptions as in \Prop{} 1.1. The
\res{}s of the
\op{} P in (1.1) in the rectangle $]E_0-\epsilon _0,E_0+\epsilon
_0[-i[0,\epsilon _1[$, for
$\epsilon _0,\epsilon _1>1$ \sufly{} small, are simple and labelled by
the two quantum numbers $k\in{\bf Z}$, $\ell\in{\bf N}$, and they are of
the form 
\eekv{2.37}
{
\sim \sum_{j=0}^\infty  g_j(kh-{S\over 2\pi },(\ell +{1\over 2})h)h^j
\hbox{ in the orientable case,} } {
\sim \sum_{j=0}^\infty  g_j((k+{\ell\over 2})h-{S\over 2\pi },(\ell
+{1\over 2})h)h^j\hbox{ in the non-orientable case,} }

Here $g_0(\tau ,\sigma )=f(\tau )-i\mu (\tau )\sigma +{\cal O}(\sigma
^2)$, where $f,\mu $ are real, $f'(\tau )>0$, $\mu >0$. Recall
that we are in the orientable case when the two \ev{}s of the
Poincar\'e map of $\gamma (E_0)$ are positive and that we are in the
non-orientable case when they are negative.
\rm\medskip

\bigskip

\centerline{\bf 3. Saddle point resonances.}
\medskip
\par Let 
\ekv{3.1}{P=-h^2\Delta +V(x),\ p(x,\xi )=\xi ^2+V(x),\ x,\xi \in {\bf
R}^2,}
where $V$ is analytic with 
\ekv{3.2}{
V(0)=E_0,\, V'(0)=0,\, {\rm sgn\,}V''(0)=(1,1),
}
so that $V''(0)$ is \nondeg{} and has one \ev{} of each sign. Assume that
the general assumptions of [HeSj] are fulfilled so that we can define the
\res{}s in a fixed \neigh{} of $E_0$, when $h>0$ is small enough. Also
assume that the union of trapped trajecories in $p^{-1}(E_0)$ is just
the point $(0,0)$. Under these assumptions a
result of [KaKe] gives all \res{}s in $D(E_0,h^\delta )$ for any fixed
$\delta >0$. In [MeSj], section 7, we got all \res{}s in a
disc $D(E_0,r_0)$ for some small but fixed $r_0$, outside small conical
\neigh{}s of $]0,\infty [$ and $-i]0,\infty [$. In this section we show
how to get the \res{}s also in such \neigh{}s.

\par After a linear change of $x$-coordinates (and the corresponding dual
change in $\xi $), we may assume that 
\ekv{3.3}
{
p(x,\xi )=E_0+p_0(x,\xi )+p_1(x,\xi )+...
}
near $(x,\xi )=(0,0)$, where $p_j(x,\xi )$ is a \hm{} \pol{} of degree
$2+j$ and 
\ekv{3.4}
{
p_0(x,\xi )={\lambda _1\over 2}(\xi _1^2-x_1^2)+{\lambda _2\over 2}(\xi
_2^2+x_2^2),\ \lambda _j>0. }
(Actually, $p_1,p_2,...$ are \indep{} of $\xi $.)

\par [KaKe] showed how to adapt the [HeSj]-theory and realize $P$ as
acting in $H(\Lambda )$-spaces, where $\Lambda \subset {\bf C}^4$ is an
IR-\mfld{} which coincides with $T^*(e^{i\pi /4}{\bf R}_{x_1}\oplus {\bf
R}_{x_2})$ near $(0,0)$ and has the property that $\forall \epsilon
>0,\exists \delta >0$ such that $(x,\xi )\in \Lambda ,\, {\rm
dist\,}((x,\xi ),(0,0))>\epsilon \Rightarrow \vert p(x,\xi )-E_0\vert
>\delta $. 

\par This means essentially (modulo an argument using a Grushin reduction
as in [MeSj]) that the study of the \res{}s of $P$ near $E_0$ can be
viewed as an \ev{} problem for $P$ after the complex scaling $x_1=e^{i\pi
/4}\widetilde{x}_1$, $\widetilde{x}_1\in{\bf R}$. 

\par The principal symbol of the scaled \op{} becomes after writing
$x_1$ instead of $\widetilde{x}_1$: 
\ekv{3.5}
{
\widetilde{p}(x,\xi )=p(e^{i\pi /4}x_1,x_2,e^{-i\pi /4}\xi _1,\xi
_2)=E_0+p_0+p_1+..., }
where the new $p_0$ is given by 
\ekv{3.6}
{p_0(x,\xi )={\lambda _1\over 2i}(\xi _1^2+x_1^2)+{\lambda _2\over 2}(\xi
_2^2+x_2^2)}
and can be identified with the restriction of the old $p_0$ to $\Lambda
$. $p_0$ takes its values in the quarter plane $\Re E\ge 0,\, \Im E\le
0$. $p_0^{-1}(E)\cap {\bf R}^4$ is a Lagrangian torus when $\Re E>0$,
$\Im E<0$ and degenerates into a closed $H_{p_0}$-trajectory, when $\Re
E$ or $\Im E$ vanishes. The contributions from these degenerate regions
were precisely the ones we did not study in [MeSj].

\par More explicitly, we have the closed trajectory 
$$\gamma_0(r_2\lambda _2):\ {1\over 2}(x_2^2+\xi _2^2)=r_2,\ x_1=\xi
_1=0,$$
for $r_2>0$, of period $2\pi /\lambda _2$ and energy $r_2\lambda _2$, and
the closed trajectory 
$$\gamma _0({r_1\lambda _1\over i}):\ {1\over 2}(x_1^2+\xi _1^2)=r_1,\
x_2=\xi _2=0,$$
for $r_1>0$, of period $2\pi i/\lambda _1$ and energy $r_1\lambda _1/i$.
We consider $\gamma _0(...)$ as a real curve with the time parameter
varying on the segment $[0,T(\gamma _0(...))]$, where $T(\gamma _0(...))$
is the (complex) period of $\gamma _0(...)$. The corresponding Poincar\'e
maps are \nondeg{} and of hyperbolic type. This implies that for $E\in{\bf
C}$, close to $r_2\lambda _2$ or $r_1\lambda _1/i$, the energy surface
$p_0^{-1}(E)$ contains a unique closed trajectory $\gamma _0(E)$ close to
$\gamma _0(r_2\lambda _2)$ or $\gamma _0(r_1\lambda _1/i)$ and of period
close to $2\pi /\lambda _2$ or $r_1\lambda _1/i$, but even when $E$
belongs to the "allowed" quadrant $\Re E>0,\, \Im E<0$, $\gamma _0(E)$
will in general not be real. Indeed, the action of $\gamma (E)$ will be
real precisely when $E=r_2\lambda _2$ or $E=r_1\lambda _1/i$ for $r_j>0$.
(Recall that the derivative of the action \wrt{} the energy is equal to
the period.)

\par Consider next the full symbol $p=\widetilde{p}$ in (3.5), for
energies $E$ with $\vert E-E_0\vert \sim \epsilon $, $0<\epsilon \ll 1$.
After the change of variables $(x,\xi )=\epsilon
^{1/2}(\widetilde{x},\widetilde{\xi })$, we get (dropping the tildes):
\ekv{3.7}
{
{p(\epsilon ^{1/2}(x,\xi ))-E_0\over \epsilon }=p_0+\epsilon ^{1/2}p_1+...
=p_0+{\cal O}(\epsilon ^{1/2})=:p(x,\xi ,\epsilon ), }
to be considered in a region with $\vert (x,\xi )\vert \sim 1$. With
$E=E_0+\epsilon F$, $\vert F\vert \sim 1$, we have that
$p^{-1}(E)$ corresponds to $p(\cdot ,\epsilon )^{-1}(F)$. Now for $F$ 
close to $r_2\lambda _2$ or to
$r_1\lambda _1/i$, we see that $p(\cdot ,\epsilon )^{-1}(F)$ contains a
closed trajectory $\gamma (F)=\gamma (F,\epsilon )$, close to $\gamma_0
(r_2\lambda _2)$ or to $\gamma _0(r_1\lambda _1/i)$. $\gamma (F)$ is an
$\epsilon ^{1/2}$-\pert{} of $\gamma _0(F)$, and we have two curves
$c_j$, $j=1,2$, in the $F$-plane at distance ${\cal O}(\epsilon ^{1/2})$
from
${\bf R}_+$ and ${1\over i}{\bf R}_+$ along which the action of $\gamma (F)$ is
real. Let $\Gamma _2,\Gamma _1$ be the corresponding unions of $\gamma
(F)$'s. Then $\Gamma _j$ are of real dimension 2.\medskip

\par\noindent \bf Lemma 3.1. \it We can find smooth IR-\mfld{}s $\Lambda
_1,\Lambda _2$ which are $\epsilon ^{1/2}$-\pert{}s of $\Lambda $, with
$\Gamma _j\subset \Lambda _j$.\rm\medskip

\par\noindent \bf Proof. \rm Fix an index $j=1$ or $2$. It is easy to see
that $\Gamma =\Gamma _j$ is totally real and that the corresponding
complexification is $\widetilde{\Gamma }=\bigcup_{E\in {\rm
neigh\,}(c_j,{\bf C})}\gamma (E)$.
$\widetilde{\Gamma }$ is a symplectic \mfld{} with complex symplectic
coordinates given by $t(\rho ),p(\rho )$, for $\rho \in \widetilde{\Gamma
}$, where $t(\rho )$ is defined by $\rho =\exp (t(\rho )H_p)(\omega (\rho
))$, where $\omega (\rho )$ belongs to some 1-dimensional "initial"
\mfld{} $W\subset \widetilde{\Gamma }$ which is transversal to the $H_p$
direction. (In order to fix the ideas, we take $W$ to be the intersection
of $\widetilde{\Gamma }$ with $x_1=\xi _1=\xi _2=0$, $\Re x_2>0$ for
$j=2$,  and with $\xi _1=x_2=\xi _2=0$, $\Re x_1 >0$, for $j=1$.)  Notice
that
$t$ will be multivalued. Thus ${\sigma _\vert}_{\widetilde{\Gamma
}}=dp\wedge dt$, if $\sigma $ denotes the complex symplectic form. On
$\widetilde{\Gamma }$, the action $S$ and period $T$ can be viewed as
\fu{}s of $p$, and a second set of symplectic coordinates on
$\widetilde{\Gamma }$ is ${t\over T(p)},S(p)$. Indeed,
$$dS(p)\wedge d{t\over T(p)}=S'(p)dp\wedge {dt\over T(p)}-{tS'(p)dp\over
T(p)^2}\wedge dT(p)=dp\wedge dt,$$
since $S'(p)=T(p)$. On $\Gamma $ both $t/T(p)$ and $S(p)$ are real, so
we see that $\Gamma $ is a \it real \rm symplectic \mfld{}. 

\par After applying a complex \canform{}, we can assume that $\Lambda $
is given by $\xi ={2\over i}{\partial \Phi _0\over \partial x}$,
$x\in{\bf C}^2$, where $\Phi _0$ is a real \stpsh{} quadratic form. (We
only consider the part of $\Lambda $ which is linear.) Then $\pi
_x(\Gamma )$ is of real dimension 2 and topologically this set is a
circle. $\Gamma $ is of the form $\xi =G(x)$, $x\in \pi _x(\Gamma )$, where
$G={2\over i}{\partial \Phi _0\over \partial x}(x)+{\cal O}(\epsilon
^{1/2})$. We look for a real and smooth \fu{} $\Phi (x)$, such that
$G(x)={2\over i}{\partial \Phi \over \partial x}$ on $\pi _x(\Gamma )$. If
we had such a \fu{}, then on $\pi _x(\Gamma )$, we would get
$${\partial \Phi \over \partial x}(x)={i\over 2}G(x),\ {\partial \Phi
\over \partial \overline{x}}=-{i\over 2}\overline{G(x)},\ d\Phi (x)={i\over
2}(G(x)dx-\overline{G(x)}\overline{dx}),$$
and ${d\Phi _\vert}_{\Gamma }\simeq{i\over 2}{(\xi dx-\overline{\xi
}\overline{dx})_\vert }_\Gamma =-{\Im \xi dx_\vert}_{\Gamma }$. The
differential of the last expression is $0$, since $\Gamma $ is real
symplectic. Hence we can find $\Phi $ locally. Let $\pi _x(\gamma )$ be
the projection of one of the closed orbits, $\gamma $, that constitute
$\Gamma $. Then for the  possibly multivalued \fu{} $\Phi $, we have 
$\int_{\pi _x(\gamma )}d\Phi =-\Im \int_\gamma  \xi dx=0$, since the
actions are real. Hence we can find $\Phi $ globally and the lemma
follows.\hfill{$\#$}\medskip

\par $\Gamma _2$ will be real in the original coordinates and even a
union of hyperbolic trajectories. This fact will not be used explicitly
since we want to treat the case of
$\Gamma _1$ at the same time. Fix $j=1$ for a maximum of generality. Since 
$$p(x,\xi ,\epsilon )={\lambda _1\over 2i}(x_1^2+\xi _1^2)+{\lambda
_2\over 2}(x_2^2+\xi _2^2)+{\cal O}(\epsilon ^{1/2}),$$
we see that the complexification $\widetilde{\Gamma }_1$ is of the form
$$(x_2,\xi _2)=f_\epsilon (x_1,\xi _1)={\cal O}(\epsilon ^{1/2}),\,\,
\vert x_1\vert +\vert \xi _1\vert \sim 1.$$
The corresponding real $\Gamma _1$ is an $\epsilon ^{1/2}$-perturbation
of (a \neigh{} of $S^1$ in) ${\bf R}^2_{x_1,\xi _1}$. Recall that $\Gamma
_1$ is a union of trajectories of close to imaginary periods. Write 
$$p={1\over i}({\lambda _1\over 2}(x_1^2+\xi _1^2)-{\lambda _2\over
2i}(x_2^2+\xi _2^2))+{\cal O}(\epsilon ^{1/2}).$$
After a complex \canform{} in $x_2,\xi _2$, we can write 
$$p={1\over i}({\lambda _1\over 2}(x_1^2+\xi _1^2)-\lambda _2x_2\xi
_2)+{\cal O}(\epsilon ^{1/2}),$$
where $\Lambda $ now is given by: $x_1,\xi _1\in{\bf R}$, $\xi _2={1\over
i}\overline{x_2}$. 

\par Let $\widetilde{\Gamma }_+$, $\widetilde{\Gamma }_-$ be the complex
incoming and outgoing hypersurfaces for the flow of $iH_p$, that contain
$\widetilde{\Gamma }$. Then $\xi _2={\cal O}(\epsilon ^{1/2})$ on $\Gamma
_+$, and $x_2={\cal O}(\epsilon ^{1/2})$ on $\widetilde{\Gamma }_-$. We
can choose symplectic coordinates as in section 1 (but now in the complex
domain). First choose $t,\tau $ on $\widetilde{\Gamma }$ with $t,\tau $
real on $\Gamma $, such that ${\sigma _\vert}_{\widetilde{\Gamma }}=d\tau
\wedge dt$, ${\rm ext\,}t=t+2\pi $ and $\tau =\tau (p)$.  (They are
essentially the action angle coordinates "$t/T,S$" in the proof of the
Lemma.) Then choose a \hol{} \fu{} $\xi $ with $\xi =\xi _2+{\cal
O}(\epsilon ^{1/2})$, such that ${\xi _\vert}_{\Gamma _+}=0$. (We are now
in the orientable case, $\xi $ will be single-valued.) Then, let
$x=x_2+{\cal O}(\epsilon ^{1/2})$ solve $H_\xi x=1$, ${x_\vert}_{\Gamma
_-}=0$, and finally extend $t,\tau $ to a full \neigh{} of
$\widetilde{\Gamma }$ as solutions of 
$$H_x t =H_\xi t=0,\ H_x\tau =H_\xi \tau =0.$$

\par Then $(t,\tau ,x,\xi )$ are symplectic coordinates, and 
\ekv{3.8}
{ip(t,\tau ,x,\xi ,\epsilon )=f_\epsilon (\tau )-\widetilde{\mu
}_\epsilon (t,\tau ,x,\xi )x\xi ,}
where $f_\epsilon =f_0+{\cal O}(\epsilon ^{1/2})$, $\widetilde{\mu }
_\epsilon =\mu _0(\tau ,x,\xi )+{\cal O}(\epsilon ^{1/2})$, and
$f_0,\mu_0 $ have the same properties as $f,\mu $ in section 1. As in
that section, we can make $\widetilde{\mu }_\epsilon $ \indep{} of $t$
(up to \ably{} high order in $(x,\xi )$), after composition with a
\canform{} close to the identity. The IR-\mfld{} of Lemma 3.1 can be taken
to be 
\ekv{3.9}{t\in S^1,\, \tau \in ]{1\over 2},{3\over 2}[,\, \xi ={1\over
i}\overline{x}, x\in{\bf C},\, \vert x\vert \le 1/C.}

\par In order to apply \Th{} 2.3, we consider $\epsilon
^{-1}(P(x,hD_x;h)-E_0)$ with leading symbol (3.7) and write $x=\epsilon
^{1/2}\widetilde{x},\, hD_x=\widetilde{h}D_{\widetilde{x}},\,
\widetilde{h}=h/\epsilon $. Then 
\ekv{3.10}
{
{1\over \epsilon
}(P(x,hD_x;h)-E_0)=\widetilde{P}(\widetilde{x},\widetilde{h}D_{\widetilde{x}},
\epsilon ;h). }
It is now clear that the conclusion of \Th{} 2.3 applies and we get (in
the case $j=1$):
\medskip
\par\noindent \bf \Prop{} 3.2. \it For $\epsilon $ and $h/\epsilon $
small enough, the \res{}s of $P$ in a rectangle $]E_0-{\epsilon \over
C},E_0+{\epsilon \over C}[-i]{\epsilon \over 2},2\epsilon [$ are simple
 and of the
form 
\ekv{3.11}{
\sim E_0+\sum_{j=0}^\infty  p_j(k{h\over \epsilon }-{1\over 2},(\ell
+{1\over 2}){h\over \epsilon },\epsilon )({h\over \epsilon })^j, }
labelled by two quantum numbers $k\in{\bf Z},\ell\in{\bf N}$, where
$p_j$ is smooth in $\tau ,\sigma ,\epsilon ^{1/2}$ in a \neigh{} of
$(0,0,0)\in{\bf R}^2\times [0,\infty [$, 
$$p_0(\tau ,\sigma ,\epsilon )=\epsilon (-if(\tau )+\mu (\tau )\sigma
+{\cal O}(\epsilon ^{1/2}+\sigma ^2)).$$
Here $f,\mu $
are real, $f'(\tau )>0$, $\mu >0$.\medskip\rm

\par Here we have chosen $\tau =0$ to correspond to the closed trajectory
of action $\pi $. Replace $p_j(\tau ,\sigma ,\epsilon )$ by $p_j(\tau
+{1\over 2},\sigma ,\epsilon )$, so that (3.11) becomes
$$\sim E_0+ \sum_{j=0}^\infty  p_j(k{h\over \epsilon },(\ell +{1\over
2}){h\over
\epsilon },\epsilon )({h\over \epsilon })^j,$$
and after modification of $p_j$, $j\ge 1$:
\ekv{3.12}
{\sim E_0+ \sum_{j=0}^\infty  p_j((k+{1\over 2}){h\over \epsilon },(\ell
+{1\over 2}){h\over
\epsilon },\epsilon )({h\over \epsilon })^j.}
Use again that the \ev{}s are \indep{} of $\mu \sim 1$, if we replace
$\epsilon$ by $\mu \epsilon $ and conclude:
\ekv{3.13}
{
p_j({\tau \over \mu},{\sigma \over \mu},\mu\epsilon )\mu
^{-j}=p_j(\tau ,\sigma ,\epsilon ). }
Using this, we define
\ekv{3.14}
{
p_j(\tau ,\sigma ,1)=p_j({\tau \over \epsilon },{\sigma \over \epsilon
},\epsilon )\epsilon ^{-j}, }
for $\vert (\tau ,\sigma )\vert \sim \epsilon  $ and $(\tau ,\sigma )$ in
the rectangle of the \prop{}. Then (3.12) becomes
\ekv{3.15}
{\sim E_0\sum_{j=0}^\infty p_j((k+{1\over 2})h,(\ell +{1\over 2})h,1)h^j.}

\par We finally connect this to the results of [KaKe]. Start again with the
\op{} $P$ after the complex scaling $x_1=e^{i\pi /4}\widetilde{x}_1$,
$\widetilde{x}_1,x_2\in {\bf R}$. We have the qBnf (with $\simeq$
indicating equivalence by conjugation by an elliptic \fop{})
\ekv{3.16}
{P\simeq P^{(N)}({1\over 2}(x_1^2+(hD_{x_1})^2),{1\over
2}(x_2^2+(hD_{x_2})^2);h)+R_{N+1}(x,hD_x;h),}
where $R_{N+1}(x,\xi ;h)={\cal O}((h,x,\xi )^{N+1})$.

\par Recall that we have the \res{}s of $P$:
\ekv{3.17}
{E_0+\sum_{j=0}^\infty  p_j((k+{1\over 2}){h\over \epsilon },(\ell+{1\over
2}){h\over \epsilon },\epsilon )({h\over \epsilon })^j.}
From (3.16) and the method of obtaining (3.17), we see that
$p_j=p_j^{(2N)}+{\cal O}(\epsilon ^{N+1})$, $j\le N$, where 
\ekv{3.18}
{E_0+\sum_{j=0}^\infty  p_j^{(2N)}((k+{1\over 2}){h\over \epsilon },(\ell+{1\over
2}){h\over \epsilon },\epsilon )({h\over \epsilon })^j.}
are the corresponding \ev{}s of $P^{(2N)}$. The $p_j^{(2N)}$ also satisfy
(3.13).

\par On the other hand, we can apply the form of $P^{(2N)}$ more directly,
to see that $P^{(2N)}$ has the \ev{}s 
$$E_0+\sum_{j=0}^\infty  p_j^{(2N)}((k+{1\over 2})h,(\ell +{1\over
2})h,1) h^j,$$
where 
$$P^{(2N)}(\tau ,\sigma ;h)\sim \sum_0^\infty  p_j^{(2N)}(\tau ,\sigma
,1)h^j,$$
so $p_j^{(2N)}(\tau ,\sigma ,1)$ is smooth in $D(0,\epsilon _0)\cap
[0,\infty [^2$. For $\epsilon \sim \vert (\tau ,\sigma )\vert $ and
$j\le N$, we have 
$$p_j(\tau ,\sigma ,1)=p_j({\tau \over \epsilon },{\sigma \over \epsilon
},\epsilon )=p_j^{(2N)}({\tau \over \epsilon },{\sigma \over \epsilon
},\epsilon )+{\cal O}(\epsilon ^{N+1})=p_j^{(2N)}(\tau ,\sigma ,1)+{\cal
O}((\tau ,\sigma )^{N+1}),$$
and using also that similar relations hold for the derivatives, we see that
$p_j(\tau ,\sigma ,1)$ are smooth in $D(0,\epsilon _0)\cap [0,\infty [^2$.
\medskip

\par\noindent \bf \Th{} 3.3. \it The \res{}s of $P$ in $D(E_0,\epsilon
_0)$ are simple, labelled by $k,\ell\in{\bf N}$, and of the form
\ekv{3.19}
{
E_0+\sum_{j=0}^\infty  p_j((k+{1\over 2})h,(\ell +{1\over 2})h,1)h^j,
}
where $p_j(\tau ,\sigma ,1)\in C^\infty (D(0,\epsilon _0)\cap [0,\infty
[^2)$.\rm

\bigskip

\centerline{\bf References.}
\medskip
\par\noindent [GeSj] C. G\'erard, J. Sj\"ostrand, \it Semiclassical
resonances generated by a closed   trajectory of hyperbolic type, \rm 
Comm. Math. Phys.,108(1987), 391-421.
\smallskip
\par\noindent [HeSj] B. Helffer, J. Sj\"ostrand, \it R\'esonances en
limite semiclassique, \rm Bull. de la SMF 114(3), M\'emoire 24/25(1986)
\smallskip
\par\noindent
[Ia] A. Iantchenko, \it La forme normale de Birkhoff pour
un op\'erateur int\'egral de Fourier, \rm Asymptotic Analysis,
17(1)(1998), 71--92. 
\smallskip
\par\noindent [IaSj] A. Iantchenko, J. Sj\"ostrand, \it Birkhoff normal
forms for Fourier integral operators II, \rm Preprint Nov 2001,
http://xxx.lanl.gov/ps/math.SP/0111134, Amer. J. Math., to appear.
\smallskip
\par\noindent [KaKe] N. Kaidi, Ph. Kerdelhu\'e, \it Forme normale de
Birkhoff et r\'esonances, \rm Asympt. Anal. 23(2000), 1--21.
\smallskip
\par\noindent [MeSj] A. Melin, J. Sj\"ostrand, \it Bohr-Sommerfeld
quantization condition for non-\sa{} \op{}s in dimension 2. \rm Preprint
juin 2001, http://xxx.lanl.gov/ps/math.SP/0111293
\smallskip
\par\noindent [Sj] J. Sj\"ostrand, \it Semiclassical resonances generated
by a non-degenerate critical   point, \rm Springer LNM, 1256, 402-429.
\smallskip
\par\noindent [Sj2] J. Sj\"ostrand, \it Geometric
bounds on the density of resonances for semiclassical problems, \rm Duke
Mathematical Journal, 60(1)(1990), 1-57.
\smallskip
\par\noindent [Sj3] J. Sj\"ostrand, \it Singularit\'es analytiques
microlocales, \rm Ast\'erisque, 95(1982).
\smallskip 
\par\noindent [Sj4] J. Sj\"ostrand, \it Semi-excited states in
non-degenerate potential wells, \rm Asymptotic Analysis, 6(1992), 29-43.
\smallskip
\par\noindent [SjZw] J. Sj\"ostrand, M. Zworski, \it Quantum monodromy and
semi-classical trace formulae, \rm J. Math. Pures
Appl., 81(2002), 1--33.
\end